\documentclass[journal]{IEEEtran}

\newcommand{\superscript}[1]{\ensuremath{^{\textrm{#1}}}}

\usepackage{multicol}
\usepackage{setspace}
\onecolumn

\usepackage{cite}

\ifCLASSINFOpdf
   \usepackage[pdftex]{graphicx}
\else
\fi
\usepackage[cmex10]{amsmath}
\usepackage{algorithm}
\usepackage{algpseudocode}

\usepackage{array}

\ifCLASSOPTIONcompsoc
  \usepackage[caption=false,font=normalsize,labelfont=sf,textfont=sf]{subfig}
\else
  \usepackage[caption=false,font=footnotesize]{subfig}
\fi
\usepackage {siunitx}
\usepackage{ifthen}
\usepackage{calc}
\usepackage{makeidx}
\usepackage{nomencl}

\newcommand{\meru}[1]{}

\makenomenclature

\nomenclature[aa]{$n\in N$}{Supplier $n$ and set of $N$ suppliers.}
\nomenclature[ab]{$m\in M$}{Load aggregator $m$ and set of $M$ load aggregators.}
\nomenclature[ac]{$h\in H$}{Time slot $h$ and set of $H$ slots in the horizon.}
\nomenclature[ad]{$r_i\in R$}{Node $i$ in the system and set of $R$ Nodes.}
\nomenclature[ae]{$i,j$}{Indices.}   
\nomenclature[af]{$\alpha,\beta$}{Incremental changes for $P_\text{B,cha}^{\text{max},m}$ and $E_\text{loss}^{\text{B},m}$ calculation.} 
\nomenclature[af1]{$\Delta h$}{Time step interval.} 

\nomenclature[ag]{$P_\text{G}^n(h)$}{Generated active power by supplier $n$ in slot $h$.}
\nomenclature[ah]{$C_\text{G}^n(P_\text{G}^n)$}{Cost function of supplier $n$.}
\nomenclature[ai]{$\rho_\text{G}^n(P_\text{G}^n)$}{Bid of supplier $n$ for generating $P_\text{G}^n$.}
\nomenclature[aj]{$P_\text{G}^{\text{max},n}$}{Maximum power limit of supplier $n$.}
\nomenclature[ak]{$P_\text{G}^{\text{min},n}$}{Minimum power limit of supplier $n$.}

\nomenclature[al]{$E^m$}{Total energy requirement of aggregator $m$ over a horizon.}
\nomenclature[am]{$P_\text{F}^m(h)$}{Aggregated flexible power demand of aggregator $m$ in slot $h$.}
\nomenclature[an]{$P_\text{F}^{\text{max},m}$}{Maximum flexible demand limit of aggregator $m$.}
\nomenclature[ao]{$P_\text{F}^{\text{min},m}$}{Minimum flexible demand limit of aggregator $m$.}
\nomenclature[ao1]{$P_\text{L}^m(h)$}{Aggregated inflexible demand of aggregator $m$ in slot $h$.}
\nomenclature[ao2]{$P_\text{U}^m(h)$}{Aggregated demand of price-responsive users communicating with aggregator $m$ before reshaping in slot $h$.}
\nomenclature[ao3]{$P_\text{PV}^m(h)$}{Aggregated PV generation of price-responsive users communicating with aggregator $m$ in slot $h$.}
\nomenclature[ap]{$P_\text{LF}^m(h)$}{Aggregated nett demand of aggregator $m$ in slot $h$.}

\nomenclature[aq]{$P_\text{B}^{m}(h)$}{Battery storage power for aggregator $m$ in slot $h$.}
\nomenclature[aq1]{$P_\text{B,cha}^{\text{max},m}$}{Maximum battery charging rate for aggregator $m$. This is a limiting variable to ensure that the total storage capacity is not exceeded. }
\nomenclature[ar]{$B_\text{SOC}^{m}(h)$}{Battery storage state of charge (SOC) of price-responsive users communicating with aggregator $m$.}
\nomenclature[as]{$B^{\text{max},m}_\text{SOC}$}{Maximum battery storage SOC limit of price-responsive users communicating with aggregator $m$.}
\nomenclature[at]{$B^{\text{min},m}_\text{SOC}$}{Minimum battery storage SOC limit of price-responsive users communicating with aggregator $m$.} 
\nomenclature[at1]{$E_\text{loss}^{\text{B},m}$}{Total battery energy loss of price-responsive users communicating with aggregator $m$ over a horizon.}

\nomenclature[az]{$P_\text{L}^{r_i,r_j}(h)$}{Transferred power by line from $r_i$ to $r_j$ in slot $h$.}
\nomenclature[azz]{$P_\text{L,max}^{r_i,r_j}$}{Maximum line power limit from $r_i$ to $r_j$.}
\nomenclature[azzz]{$P_\text{L,min}^{r_i,r_j}$}{Minimum line power limit from $r_i$ to $r_j$.}
\nomenclature[azzzz]{$B_{r_i,r_j}$}{Susceptance of line between  $r_i$ and $r_j$.}
\nomenclature[azzzzz]{$\delta_{r_i}(h)$}{Voltage angle at node $r_i$ in slot $h$.}

\makeatletter
\def\BState{\State\hskip-\ALG@thistlm}
\makeatother
\begin{document}
\title{Aggregated Demand Response Modelling for Future Grid Scenarios}
%
%
\author{Hesamoddin~Marzooghi,~
        Gregor~Verbi\v{c},~
        and~David~J.~Hill
\thanks{Hesamoddin~Marzooghi, Gregor~Verbi\v{c} and David~J.~Hill are with the School of Electrical and Information Engineering, The University of Sydney, Sydney, New South Wales, Australia. e-mails: ({hesamoddin.marzooghi, gregor.verbic, david.hill}@sydney.edu.au).}
\thanks{David~J.~Hill is also with the Department of Electrical and Electronic Engineering, The University of Hong Kong, Hong Kong. e-mail: (dhill@eee.hku.hk).}
\thanks{}}

\maketitle
\begin{abstract}
With the increased penetration of intermittent renewable energy sources (RESs) in future grids (FGs), balancing between supply and demand will become more dependent on demand response (DR) and energy storage. Thus, FG feasibility studies will need to consider DR for modelling nett future demand. Against this backdrop, this paper proposes a demand model which integrates the aggregated effect of DR in a simplified representation of the effect of market/dispatch processes aiming at minimising the overall cost of supplying electrical energy. The conventional demand model in the optimisation formulation is augmented by including the aggregated effect of numerous users equipped with rooftop photovoltaic (PV)-storage systems. The proposed model is suited for system studies at higher voltage levels in which users are assumed to be \textit{price anticipators}. As a case study, the effect of the demand model is studied on the load profile, balancing and loadability of the Australian National Electricity Market in 2020 with the increased penetration of RESs. The results are compared with the demand model in which users are assumed to be \textit{price takers}.
\end{abstract}


%
\IEEEpeerreviewmaketitle

\printnomenclature[1.3cm]

\section{Introduction}

\IEEEPARstart{f}{uture} grid (FG) feasibility studies have demonstrated that relying on high penetration of diverse renewable energy sources (RESs) is possible assuming enough backup generation and/or utility storage are available to keep the network in balance \cite{Energy2010, Elliston2012, Elliston2013, Budischak2013, Hart2011, Mason2010}. A preliminary study by the University of Melbourne Energy Research Institute has proposed a zero-carbon electrical grid for Australia in 2020 \cite{Energy2010}. The University of New South Wales researchers have analysed the viability of \SI{100}{\percent} RES scenarios considering a copper plate model for the Australian National Electricity Market (NEM) \cite{Elliston2012, Elliston2013}. They have suggested \SI{100}{\percent} RESs electricity in the NEM, at the current reliability standard, would be technologically feasible. Also, the least-cost mix of \SI{100}{\percent} RESs scenario has been determined for the future of the NEM. Similarly, the least-cost mix of high penetration of diverse RESs and conventional generation has been determined for the future of the PJM, California and New Zealand networks in \cite{Budischak2013, Hart2011, Mason2010}, respectively.

\par However, these studies have only focused on simple balancing by using a simplified grid model such as the copper plate model. On the other hand, the penetration of distributed generation (DG) has been increasing significantly in recent years, and greater penetration of battery storage is anticipated in power systems \cite{EPRI, RockyMountain, ATA, PV, AEMO2012A, AEMO2012, CSIRO}. In particular, global installed capacity of rooftop photovoltaic (PV) has increased from approximately 4 GW in 2003 to nearly 128 GW in 2013 mainly due to electricity price increases, government incentives and also worldwide drop of PV capital cost \cite{EPRI, RockyMountain}. In Australia, installed capacity of rooftop PV (which is mostly installed by residential and commercial customers) has grown from approximately 0.8 GW in 2011 to over 4 GW in 2014 \cite{PV}. Recent studies have demonstrated that users equipped with PV-battery systems will reach retail price parity in the foreseeable future in the USA grids and the NEM \cite{ EPRI, RockyMountain, ATA}. In light of these developments, a question arises how to model the aggregated nett demand (including DG, storage and demand response (DR)) to study FG scenarios.

\par While the effect of DR is neglected in most of the existing FG feasibility studies \cite{Energy2010, Elliston2012, Elliston2013, Budischak2013, Hart2011}, it is considered in few studies mainly through two different ways:

\textit{Implicit modelling}: DR is considered implicitly, but it is not reflected into the loads. For instance in \cite{Mason2010}, the effect of DR is considered through improving the capacity credit value for intermittent RESs (i.e. intermittency of RESs is decreased). However, due to the significant effect of loads on performance and stability of power systems, it can be expected that incorporating DR explicitly into the load models will affect the results of FG feasibility studies significantly.

\textit{Explicit modelling}: In few recent studies \cite{CSIRO, DRNEW5, Claudio, Hesam2014, DRGV, DRNEW2}, the aggregated effect of DR is reflected into the conventional demand models. Including the effect of DR into the demand models requires allowing for the interaction between demand and supply sides in some ways. This is mainly done through three different approaches. Firstly, in some studies, the supply-side is modelled physically while price-responsive users are not represented physically \cite {DRNEW5, Claudio}. In \cite{DRNEW5}, the effect of flexible loads is analysed on reserve markets. That study presumes flexible load to be represented precisely by a tank model. Also, the reserve market is too simplified and physical constraints of the electrical grids (e.g. line limits) are not considered. In \cite {Claudio}, flexible demand is represented via a price-elasticity matrix. The elasticities are a measure of the change in demand in response to a change in the electricity price, and are typically obtained from the analysis of historical data. Secondly, there are some other studies that model demand-side technologies physically while the supply-side is represented through the electricity price profile \cite{Hesam2014}. That model assumes user to be price takers, i.e. the effect of user actions is not considered in the electricity price. Thirdly, in few recent studies, both demand-side technologies and supply-side are modelled physically \cite{DRGV, DRNEW2}. This approach necessitates the need for integrated simulations in which both supply and demand sides need to be optimised jointly. As discussed in \cite{DRGV}, this approach can provide more realistic results. In \cite {DRNEW2}, the aggregated charging management approaches for plug-in electrical vehicles (PEVs) is integrated into the market clearing process. The market process, however, is too simplified and physical constraints of the electrical grids are not considered.

\par Although those models have proven their merits, a generic modelling approach is still required which can be used for any granularity level in the grid (e.g. from a city, a state or even the whole network). Furthermore, for studying FG scenarios of uptake of the various demand-side technologies, it is necessary to have a model which can easily integrate them. Against this backdrop, in this paper, we make a further step by proposing a model that integrates the aggregated effect of DR in a simplified representation of the effect of market/dispatch processes. The proposed optimisation formulation aims at minimising the overall cost of supplying electrical energy in which the conventional demand model is augmented by including the aggregated effect of numerous price anticipating users equipped with PV-storage systems. Due to the price anticipating assumption, load aggregators are considered implicitly. As a case study, the effect of the demand model on the load profile, balancing and loadability of the NEM is studied using a 14-generator model \cite{Gibbard2010}. The results are also compared with the demand model in \cite{Hesam2014} in which users are assumed to be price takers. Eight scenarios are analysed with one business as usual (BAU) and seven different DR scenarios with renewable integration. For the BAU Scenario in 2020, the electricity supply is dominated by coal, gas, hydro, and biomass; and in the Renewable Scenarios, some of the conventional coal generators in Queensland and South Australia are replaced with wind farms (WFs) and concentrated solar plants (CSPs) with storage, as suggested in \cite{Energy2010, AEMO2012A} to meet Australia's RES target \cite{DepartofResources}. Simulation results show that increasing the penetration of DR with price anticipating users, improves loadability with the increased intermittent supply penetration in the grid.

\par The remainder of the paper is organised as follows: Section II proposes the aggregated demand model considering DR. Section III describes the test-bed assumptions and modelling. Section IV describes simulation scenarios, and discusses simulation results. Finally, Section V gives conclusions.

\section{Aggregated Demand Model Considering DR}

\par Aggregated load models are commonly used in system studies to reflect the combined effect of numerous physical loads \cite{Kundur, Van}. These can be inspired by physical devices, e.g. using a large induction motor to represent all the motors connected, or by data-driven approaches. Conventional aggregate load models only account for the accumulated effect of numerous independent load changes and some relatively minor control actions. Including the effect of DR requires allowing for the interaction between demand and supply sides in some way, e.g. price signals. In this study, we propose a model that integrates the aggregated effect of DR in a simplified representation of the effect of market/dispatch processes inspired by the traditional unit commitment problem. The proposed optimisation problem aims at minimising the overall cost of supplying electrical energy in which the conventional demand model is augmented by including the aggregated effect of numerous users equipped with newer demand technologies. The demand model consists of two parts: (i) a fixed electricity demand profile (inflexible demand), and (ii) flexible demand equipped with demand-side technologies (we consider a large homogeneous population of residential and commercial PV-storage systems, but the model allows for an easy integration of other demand technologies as well). Also, the demand model considers the following assumptions:

\textit{Assumption 1}: Users are assumed to be price anticipators, i.e. the effect of users’ actions is considered on the electricity price by the load aggregators. Due to the price anticipating assumption, load aggregators are considered implicitly. We assume that users are incentivised by the load aggregators through a proper price signal derived using mechanism design \cite{Samadi2012}.

\textit{Assumption 2}: Aggregators do not aim to change the total energy consumption of the users, but instead to systematically manage and shift it. Furthermore, all load aggregators are connected to not only the power grid but also to a communication infrastructure to enable two-way communication with the users. 

\textit{Assumption 3}: Users have smart meters equipped with smart home energy management systems (SHEMSs). SHEMS implements an algorithm to schedule energy resources and storage, and so achieves DR.

\subsection{Optimisation model}

\par The optimisation model aims at minimising the overall electricity cost. The objective function of the model can be written as:
\begin{equation}
             \underset{}{\text{min}}
          \sum_{h=1}^{H}{\sum_{n=1}^{N} C_\text{G}^{n}(P_\text{G}^{n}(h))},\quad\\
\end{equation}
where, each decision horizon for the model (i.e. 24-hour period) is divided into one hour time-steps, giving a total of $H=24$ time-steps; denote a particular time-step by $h$, subject to the following constraints:

\subsubsection{Power generation limit} Generation of each supplier is a decision variable, and is constrained between the minimum and the maximum power limits as follows:
\begin{equation}	
P_\text{G}^{\text{min},n}\leq P_\text{G}^{n}(h)\leq P_\text{G}^{\text{max},n}\quad\forall h,n,\\
\end{equation}		

\subsubsection{Flexible demand, storage and PV} The following set of equations augment the conventional demand model by including the aggregated effect of numerous users equipped with PV-storage systems.
\begin{subequations}
      \begin{align}	
          P_\text{F}^{\text{min},m}\leq P_\text{F}^{m}(h)\leq P_\text{F}^{\text{max},m}\quad\forall h,m,\\
					\sum_{h=1}^{H}[P_\text{LF}^{m}(h)-P_\text{PV}^{m}(h)] \Delta h=E^m+E_\text{loss}^{\text{B},m}\quad\forall m,
		  \end{align}
\end{subequations}

\par
\begin{figure} 
\centering
\includegraphics[width=8.8cm, height=9cm]{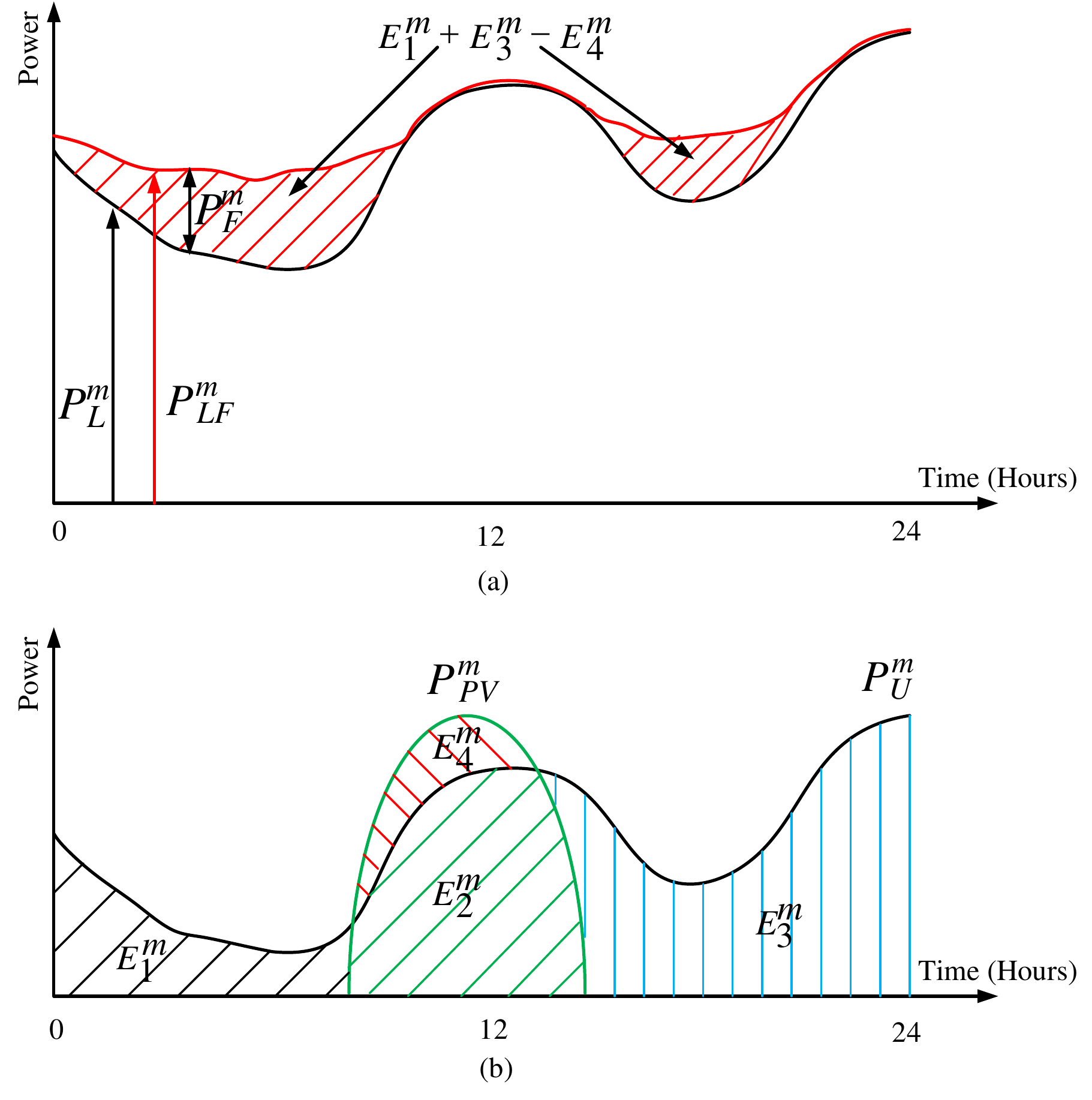} 
\caption{Simple illustration of the demand profile (a) aggregated flexible demand, inflexible demand and new nett demand profile, and (b) aggregated price-responsive demand and PV generation.}
\label{figure:LoadModel}
\end{figure}

 Flexible demand of each load aggregator, $P_\text{F}^{m}$, is a decision variable which reflects the aggregated effect of DR, and is constrained between the minimum and the maximum limits in (3a). The overall energy balance over time horizon $H$ is given by (3b). Fig.~\ref{figure:LoadModel} shows a simple illustration of the demand profile. As it can be seen in Fig.~\ref{figure:LoadModel}a, the aggregated nett demand of each load aggregator, $P_\text{LF}^{m}$, is equal to the sum of inflexible and flexible demands, i.e. $P_\text{LF}^{m}(h)=P_\text{L}^{m}(h)+P_\text{F}^{m}(h)$. When $P_\text{F}^{\text{min},m}$ is equal or more than zero ($P_\text{L}^{m}\leq P_\text{LF}^{m}$), the model represent a situation where price-responsive users would not send power back to the grid. This situation might happen in the future if feed-in tariffs (FiTs) are much less than the retail tariffs paid by the users. On the other hand, if there will be a reasonable incentive for users in the future to send power back to the grid, $P_\text{F}^{\text{min},m}$ can be relaxed to a value below zero ($P_\text{LF}^{m}<P_\text{L}^{m}$) to reflect such a situation. The flexibility of loads is due to battery storage which is modelled implicitly by considering the upper limit of the flexible load power as $P_\text{F}^{\text{max},m}(h)=P_\text{B,cha}^{\text{max},m}+P_\text{U}^m(h)$. Note that $P_\text{B,cha}^{\text{max},m}$ is a limiting variable that ensures that the total storage capacity is not exceeded, and it does not represent a physical property of a particular battery technology. This variable is determined heuristically, as explained later in this section. The battery storage power, $P_\text{B}^{m}$, is a byproduct of the optimisation problem which can be calculated as $P_\text{B}^{m}(h)= P_\text{F}^{m}(h)-P_\text{U}^m(h)-P_\text{PV}^m(h)$. Fig.~\ref{figure:LoadModel}b shows the aggregated flexible demand equipped with PV-storage systems. As it is demonstrated in Fig.~\ref{figure:LoadModel}b, $P_\text{U}^m(h)$ is the aggregated demand of price-responsive users before utilising PV-storage systems.

\par The aggregated flexible demand is determined by each aggregator in (3b) in a way that the total energy requirement for that aggregator remains constant (\textit{Assumption 2}) considering aggregated PV generation and also battery storage energy loss. Total energy requirement of aggregator $m$ over the horizon can be written as $E^m=\sum_{h=1}^{H}(P_\text{L}^m(h)+P_\text{U}^m(h))\Delta h$. $E_\text{loss}^{\text{B},m}$ in (3b) guarantees that battery round-trip efficiency is taken into account. As it is shown in Fig.~\ref{figure:LoadModel}b, part of the required energy for price responsive users is provided by PV. So, the rest of the energy has to be supplied from the grid (i.e. $E_1^m+E_3^m -E_4^m $). Price-responsive users use their battery storage to shift their consumption to cheaper time slots to minimise the overall cost of supplying electrical energy (i.e. $E_1^m +E_3^m -E_4^m $ can be spread out over the horizon due to enough storage capacity of price-responsive users in a way that the overall electricity cost minimises, the red area in Fig.~\ref{figure:LoadModel}a). In other words, the total flexible demand energy over the horizon is equal to the energy which has to be supplied from the grid plus battery storage energy loss, i.e. $\sum_{h=1}^{H}P_\text{F}^m(h)\Delta h= \sum_{h=1}^{H}(P_\text{U}^m(h)-P_\text{PV}^m(h))\Delta h+E_\text{loss}^{\text{B},m}=E_1^m+E_3^m-E_4^m+E_\text{loss}^{\text{B},m}$.

\subsubsection{Demand supply balance} Ignoring the losses in the grid, the power balancing equation can be written as (4a): 
\begin{subequations}
      \begin{align}	
          \sum_{n\in r_i}P_\text{G}^n(h)-\sum_{m\in r_i}P_\text{LF}^m(h)=\sum_{r_j}P_\text{L}^{r_i,r_j}(h)\quad\forall h, r_i\in R, \\
					P_\text{L}^{r_i,r_j}(h)=B_{r_i,r_j}(\delta_{r_i}(h)-\delta_{r_j}(h))\quad\forall h, (r_i,r_j)\in R,\\
					P_\text{L,min}^{r_i,r_j}\leq P_\text{L}^{r_i,r_j}(h)\leq P_\text{L,max}^{r_i,r_j}\quad\forall h, (r_i,r_j)\in R,
		  \end{align}
\end{subequations}
where, the power transferred between different nodes in the system is given by (4b), and is constrained by the line limits in (4c).

\par To solve the above optimisation problem, aggregated demand, PV power, maximum battery storage charging rate and corresponding battery storage energy loss are required for each load aggregator. The next subsection describes how the maximum battery storage charging rate and corresponding energy loss can be calculated for the model.

\subsection{Modelling battery storage}

\par In the proposed optimisation formulation, battery storage is modelled implicitly and its state of charge (SOC) is not a decision variable. However, it is important to consider battery SOC limits to make sure that the available storage capacity in the grid is not exceeded. $P_\text{B,cha}^{\text{max},m}$ is a limit parameter which ensures that the total storage capacity is not exceeded. It is determined using a heuristic search in Algorithm~\ref{BCM}.

\begin{algorithm}
\caption{: The proposed algorithm to determine $P_{\text{B,cha}}^{\text{max},m}$}\label{BCM}
\begin{algorithmic}[1]
\State Set $i=1,$
\State Initialize $P_{\text{B,cha},i}^{\text{max},m}=0, $
\State Run Algorithm 2,
\While{$ B_{\text{SOC}}^m $ is within limits (7)}
\State Solve the model (1-4).
\State Calculate $B_{\text{SOC}}^m$ (6).
\State $i \leftarrow i+1,$
\State $P_{\text{B,cha},i}^{\text{max},m} \leftarrow P_{\text{B,cha},(i-1)}^{\text{max},m}+\alpha, $
\State Run Algorithm 2,
\EndWhile
\State Step 3: Return $P_{\text{B,cha},(i-1)}^{\text{max},m}$.
\end{algorithmic}
\end{algorithm}

\par The battery energy loss corresponding to $P_{\text{B,cha},i}^{\text{max},m}$ is also verified based on a heuristic search in Algorithm~\ref{BCM1}. Based on the Algorithm~\ref{BCM1}, the optimisation formulation (1-4) is solved using the aggregated demand and PV power, suppliers' data, battery storage charging rate (from Algorithm~\ref{BCM}) and initial value for battery energy loss. Then, $E_{\text{loss}}^{\text{B},m}$ can be calculated based on the storage actions as follows:
\begin{equation}
        E_{\text{loss}}^{\text{B},m}=(1- \eta)\sum_{h=1}^{H} P_{\text{B}}^m(h) \quad if P_{\text{B}}^m(h)>0,
\end{equation}

\par If the difference between the battery energy loss calculated from the equation (5) and the battery energy loss from the j\superscript{th} iteration, $E_{\text{loss},j}^{\text{B},m }$, is less than the error, $\varepsilon$, search stops for that horizon. Otherwise, the battery energy loss value in the j\superscript{th} iteration changes in a small step ($\beta$) and the optimisation formulation solves again. The search continues until the difference between those two is less than the error. 

\begin{algorithm}
\caption{: The proposed algorithm to verify $E_{\text{loss}}^{\text{B},m}$ corresponding to $P_{\text{B,cha},i }^{\text{max},m }$}\label{BCM1}
\begin{algorithmic}[1]
\State Set $j=1,$
\State Initialize $E_{\text{loss},j}^{\text{B},m}=0, $
    \While{$ E_{\text{loss}}^{\text{B},m}- E_{\text{loss},j}^{\text{B},m}>\varepsilon $}
        \State Solve the model (1-4).
        \State Calculate $ E_{\text{loss}}^{\text{B},m}$ (5).
        \State $j \leftarrow j+1,$
        \State $ E_{\text{loss},j }^{\text{B},m } \leftarrow E_{\text{loss},(j -1)}^{\text{B},m }+\beta, $
  \EndWhile
 \State Step 4: Return $E_{\text{loss},j}^{\text{B},m}$.
\end{algorithmic}
\end{algorithm}

\par After battery energy loss calculation from Algorithm~\ref{BCM1}, the optimisation formulation can be solved in Algorithm~\ref{BCM} for determining $P_{\text{B,cha}}^{\text{max},m}$. To verify $P_{\text{B,cha}}^{\text{max},m}$, battery storage SOC, $B_{\text{SOC}}^m$, which is a byproduct of the optimisation problem needs to be calculated:
\begin{subequations}
      \begin{align}
       &\underset{}{}
       & & B_{\text{SOC}}^m(1) = B_{\text{SOC}}^{i,m},\quad\\
       &
       & & B_{\text{SOC}}^m(h) = B_{\text{SOC}}^m(h-1) + P_{\text{B}}^m(h-1) \quad\forall h>1,
       \end{align}
\end{subequations}        
and compared with the battery storage SOC limits: 
\begin{equation}
       B^{\text{min},m}_{\text{SOC}}\leq B_{\text{SOC}}^m(h)\leq B^{\text{max},m}_{\text{SOC}}\quad\forall h>1,
\end{equation}

\par While $B_{\text{SOC}}^m$ is within the limits, $P_{\text{B,cha}}^{\text{max},m}$ increases in a small step ($\alpha$), $E_{\text{loss}}^{\text{B},m}$ corresponding to $P_{\text{B,cha},i}^{\text{max},m}$ calculates from Algorithm~\ref{BCM1}, and the optimisation formulation solves again. This procedure repeats until $B_{\text{SOC}}^m$ violates the limit. From a step before violation, $P_{\text{B,cha}}^{\text{max},m}$ can be obtained.

\section{The Australian NEM Model}

\par
\begin{figure} 
\centering
\includegraphics[width=11cm, height=15.5cm]{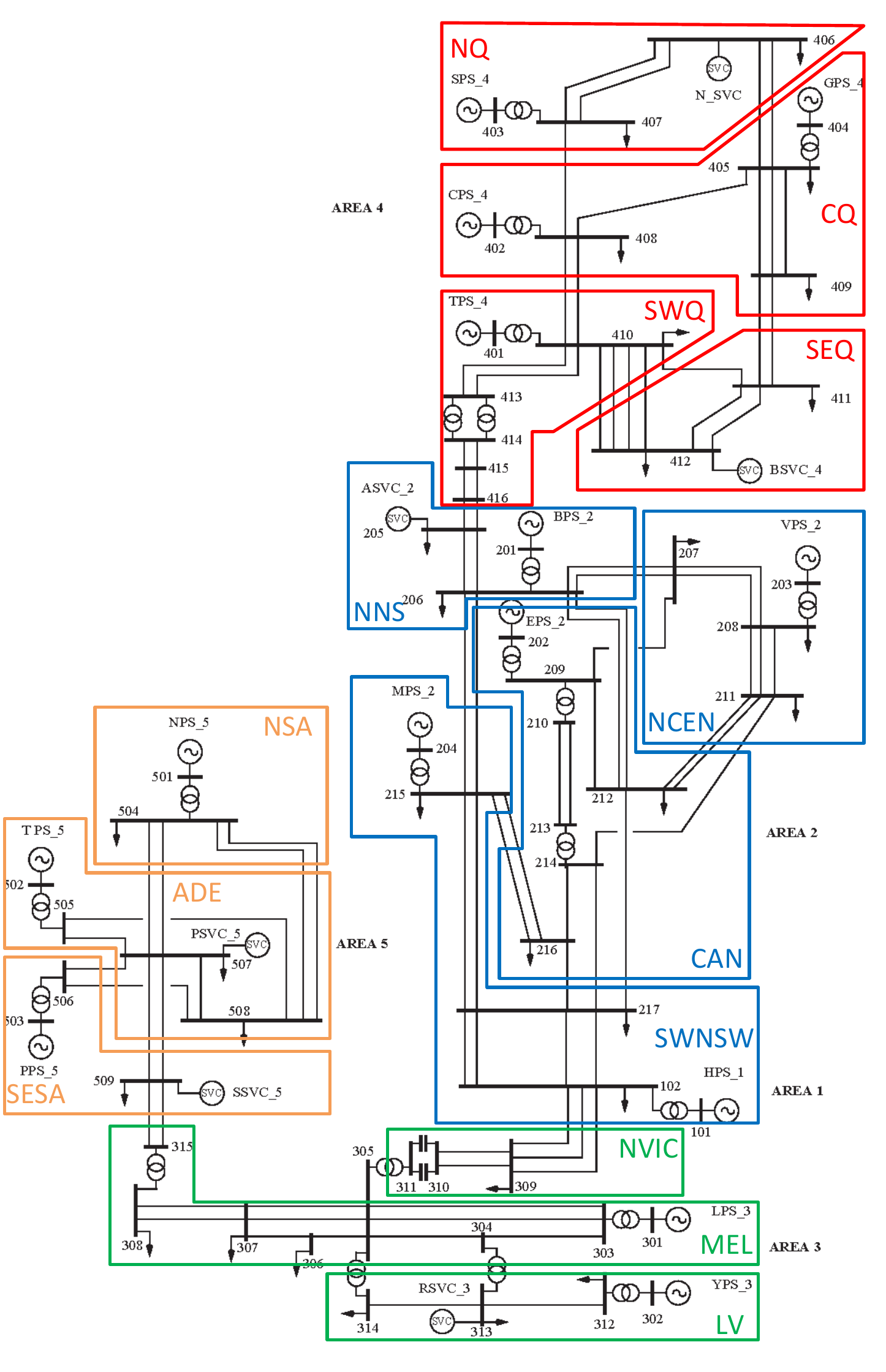}
\caption{14-generator model of the NEM.}
\label{figure:14-generator model of the NEM}
\end{figure}

 A 14-generator model of the NEM, which was originally proposed for small signal stability studies \cite{Gibbard2010}, is used as the test-bed. The schematic diagram of the 14-generator model of the NEM is shown in Fig.~\ref{figure:14-generator model of the NEM}. Areas 1 to 5 represent Snowy Hydro (SH), New South Wales (NSW), Victoria (VIC), Queensland (QLD) and South Australia (SA), respectively. The Australian Electricity Market Operator (AEMO) has proposed 16 zones for the NEM to capture differences in generation technology capabilities, costs, weather and so on in the future \cite{AEMO2012A}. In order to extract data for the proposed model and generators in 2020, the 14-generator model is matched with the 16 zones, as shown in Fig.~\ref{figure:14-generator model of the NEM}.

\par The generator technologies and the test-bed assumptions in this study follow Reference \cite{Hesam2014}. In power system studies, generators' cost functions are often modelled as quadratic. In this study, it is assumed that conventional generators submit their block bidding curves, while RESs are assumed to bid at zero cost. So, the cost function of the conventional generators can be approximated by a piece-wise linear function as:

\begin{equation}	  
     C_\text{G}^n(P_\text{G}^n(h))=\rho_\text{G}^n(P_\text{G}^n(h)) P_\text{G}^n(h) \quad\forall n,                                                                                            
\end{equation}
where, $\rho_\text{G}^n(P_\text{G}^n(h))$ are suppliers' bids. Generators' bids mimic the bidding behaviour of the generators in the NEM and use the predicted fuel price, thermal efficiency, and variable O\&M costs 2020 \cite{AEMO2012A, Energy2012}. Table~\ref{tab:The predicted cost of the chosen technologies in 2020} lists bids of the suppliers ($\rho1$, $\rho2$ and $\rho3$) in the 14-generator model in 2020. 
 
\par Each region of the 14 generator model (i.e. QLD, NSW, VIC and SA) is considered as a node for the demand model. Consequently, interstate lines are considered explicitly in the demand model. However, balancing and loadability studies are done using the 14 generator model shown in Fig.~\ref{figure:14-generator model of the NEM}. The dispatch results from the market are used for balancing and loadability studies in DIgSILENT Power Factory. In the balancing studies, if supply cannot meet the demand, the hour is recorded as unserved hour. However, if available generation exceeds demand (i.e. due to high generation of RESs), the surplus power is recorded as dumped energy and that hour is marked as a dumped hour.

\begin{table} 
\centering
\caption{the generators' bids in 2020}
\begin{tabular}{| c | c |  c | c | c |c | c |}\hline
 Gen. & Type  & AEMO   &  $\rho1$    &  $\rho2$    & $\rho3$\\
           &       & zone   & (\$/MWh)  & (\$/MWh)  & (\$/MWh)\\\hline
 BPS\_2    &  Coal     & NNS       & 28.45&  42.66    &  56.90\\
 EPS\_2    &  GT       & CAN       & 69.20&  346.0    &  692.0\\
 MPS\_2    &  Coal     & SWNSW     & 27.43&  41.15    &  54.86\\
 VPS\_2    &  Coal     & NCEN      & 26.40&  39.60    &  52.80\\
 LPS\_3    &  Biomass  & MEL       & 39.50&  59.25    &  79.00\\
 YPS\_3    &  Coal     & LV        & 21.88&  32.82    &  43.76\\
 CPS\_4    &  Coal     & CQ        & 26.14&  39.21    &  52.28\\
 GPS\_4    &  Coal     & CQ        & 26.14&  39.21    &  52.28\\
 SPS\_4    &  Coal     & NQ        & 32.74&  49.11    &  65.48\\
 TPS\_4    &  GT       & SWQ       & 73.84&  369.2    &  738.4\\
 NPS\_5    &  Coal     & NSA       & 30.89&  46.34    &  61.78\\
 PPS\_5    &  Coal     & SESA      & 30.89&  46.34    &  61.78\\
 TPS\_5    &  GT       & ADE       & 69.20&  346.0    &  692.0\\\hline         
\end{tabular}
\label{tab:The predicted cost of the chosen technologies in 2020}
\end{table}

\section{Simulation Scenarios and Results}

\par The effect of price anticipating and price taking assumptions, and different DR penetrations on the load profile, balancing and loadability of the NEM in 2020 with the increased penetration of RESs in the grid is demonstrated in this section.

\subsection{Simulation scenarios}

\par Eight scenarios are analysed with one BAU and seven different DR scenarios with renewable integration. The Renewable Scenarios are analysed with the conventional load, and different levels of DR in the proposed demand model (Section II) and the demand model in \cite{Hesam2014}. For the BAU Scenario, combinations of coal, gas, hydro, and biomass are considered for the NEM to supply the load in 2020. Then, some of the conventional coal generators in QLD and SA are replaced with CSPs together with storage and WFs, respectively to meet Australia's RES target. Displacement of conventional generators in the Renewable Scenarios and chosen capacities for RESs are inspired by studies in \cite{AEMO2012A, Energy2010}. NPS\_5 in SA is replaced with a WF with the capacity of \SI{3}{\giga\watt} using NSA data \cite{AEMO2012A}. Also, SPS\_4 and GPS\_4 in QLD are replaced with two CSPs with the capacity of \SI{4.5}{\giga\watt} each and using NQ and CQ data \cite{AEMO2012A}, respectively. It was found that delaying CSP output by 12 hours minimises the unserved and dumped energy. The RESs serve about \SI{20}{\percent} of the total demand in the Renewable Scenarios. 

\par In this study, hourly demand and PV power for the proposed demand model are obtained from the AEMO predications for 2020 \cite{AEMO2012A}. Also, DR is considered for the residential and commercial customers, i.e. \SI{60}{\percent} of the total system load in the NEM in 2020 \cite{AEMO2012A}. The industrial customers are left unaffected. Furthermore, the percentage of the residential and commercial customers with PV are considered \SI{20}{\percent}, \SI{30}{\percent} and \SI{40}{\percent} for low, medium and high uptake scenarios, respectively. Table~\ref{tab:The aggregated storage and PV capacities for each region of the NEM for different DR scenarios} shows the aggregated storage and PV capacities for each region of the NEM for different uptake scenarios. The chosen PV capacities for different uptake scenarios are inspired by the AEMO study \cite{AEMO2012}. Moreover, the chosen storage capacities roughly correspond to typical PV and storage capacity for a household in Australia. Using the algorithms in Section II, $P_{\text{B,cha}}^{\text{max},m}$ is also calculated for different DR scenarios, and is reported in Table~\ref{tab:The aggregated storage and PV capacities for each region of the NEM for different DR scenarios}.

\begin{table} 
\centering
\caption{the aggregated storage and PV capacities for each region of the NEM for different uptake scenarios}
\begin{tabular}{| c | c | c | c | c|}\hline
Region & Scenario & $B^{min,m}_{SOC}$-$B^{max,m}_{SOC}$ & $P_{B,cha}^{max,m}$ & PV capacity \\
       &          &         (GWh)                       &   (GW)              &   (GW)      \\\hline
    & Low &        0.4-4.3  &  0.44  & 1.3 \\
QLD & Medium &     0.6-6.4  &  0.63  & 1.9 \\
    & High   &     0.9-8.5  &  0.83  & 2.6 \\\hline
    & Low &        0.7-6.7  &  0.62  & 2.0 \\
NSW & Medium &     1.0-10.1 &  1.01  & 3.0 \\
    & High   &     1.4-13.5 &  1.34  & 4.1 \\\hline
    & Low &        0.5-5.0  &  0.47  & 1.5 \\
VIC & Medium &     0.8-7.5  &  0.72  & 2.3 \\
    & High   &     1.0-10.0 &  0.95  & 3.0 \\\hline
    & Low &        0.1-1.2  &  0.10  & 0.3 \\
SA  & Medium &     0.2-1.7  &  0.16  & 0.5 \\
    & High   &     0.2-2.3  &  0.22  & 0.7 \\\hline    
\end{tabular}
\label{tab:The aggregated storage and PV capacities for each region of the NEM for different DR scenarios}
\end{table}

\subsection{Accuracy of the load model}

\par
\begin{figure} 
{\subfloat[]{\includegraphics[width=8.8cm, height=6.6cm]{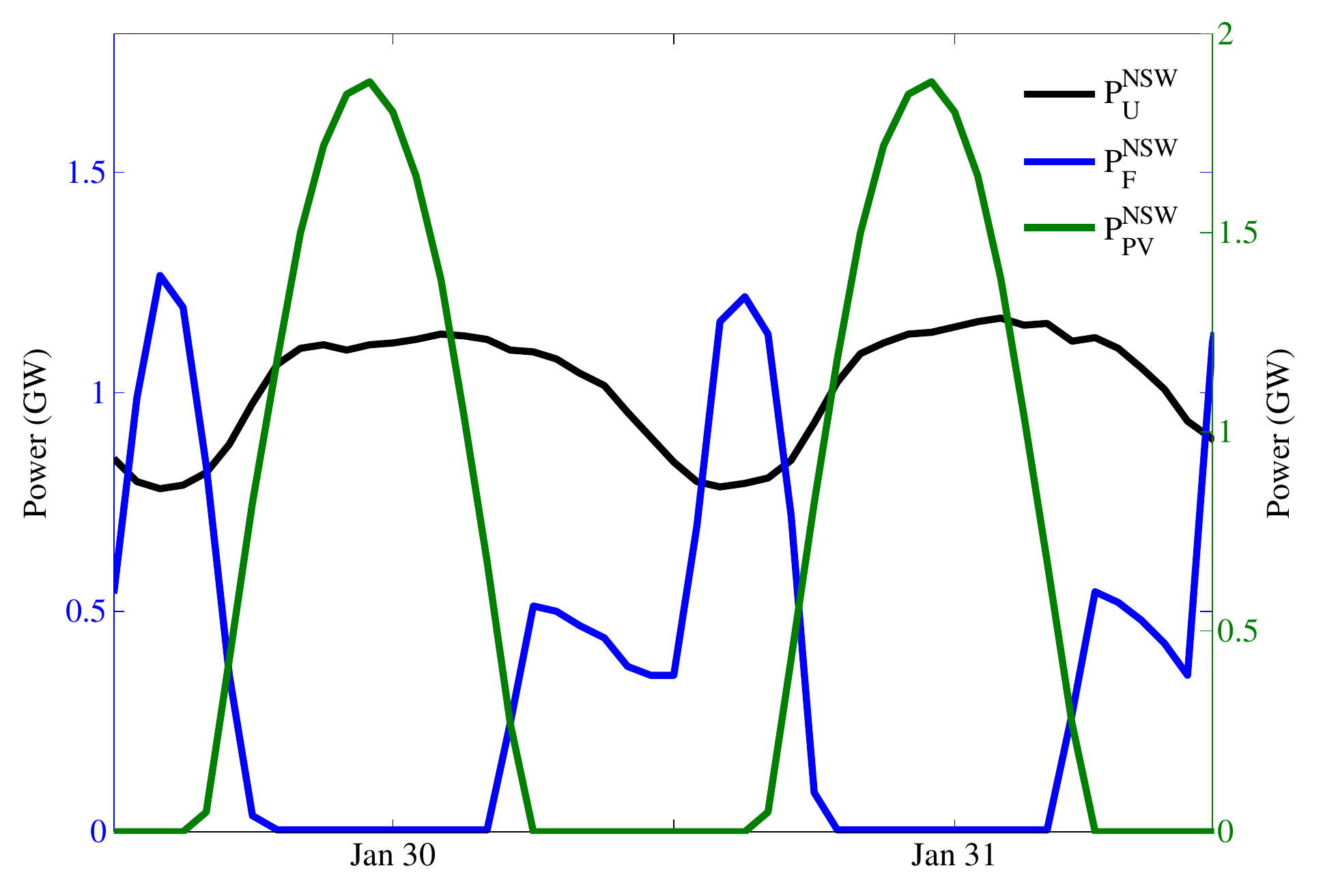}
\label{fig_second_caseA2}}} 
{\subfloat[]{\includegraphics[width=6.5cm, height=8.8cm, angle =90]{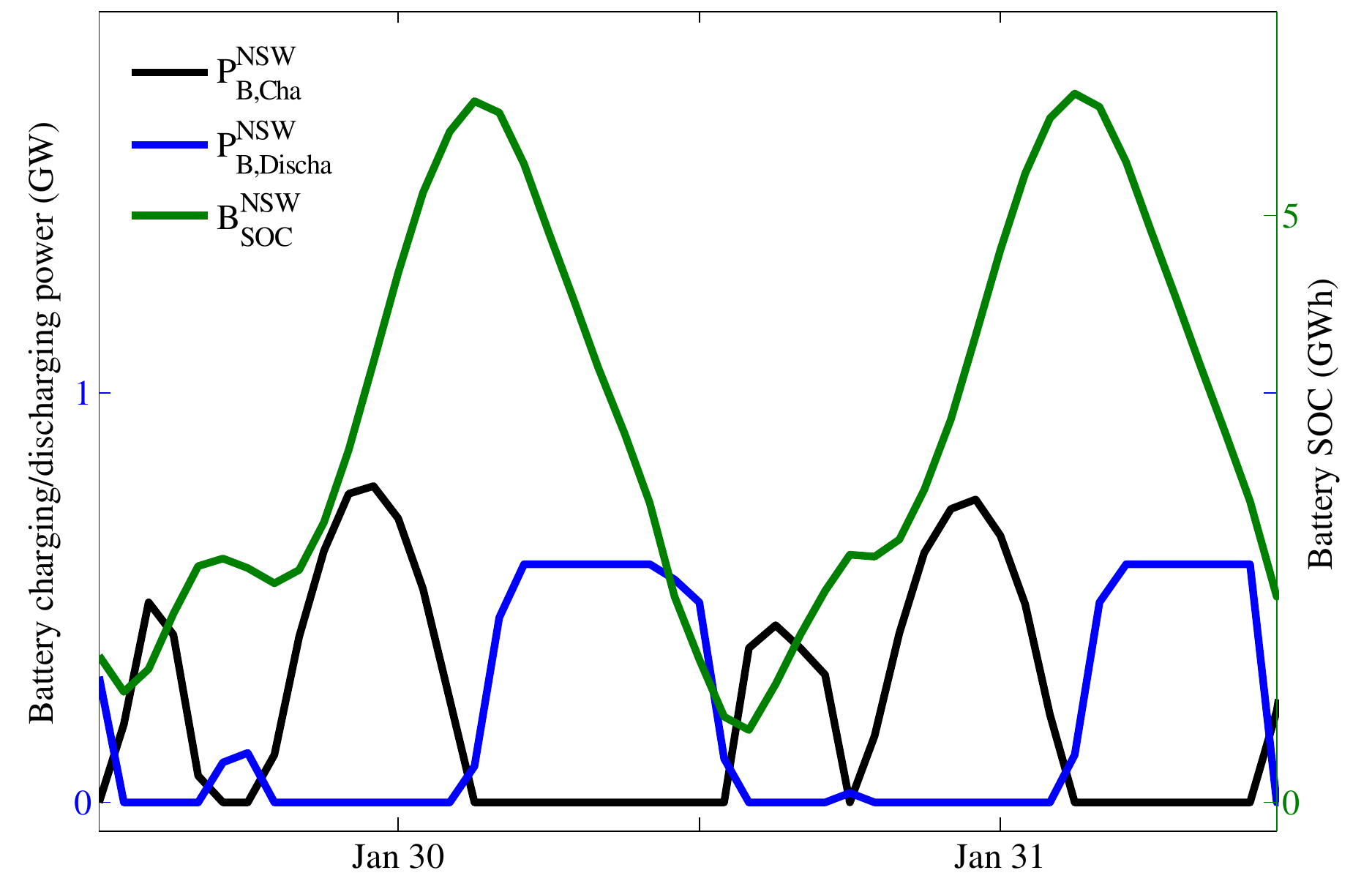}
\label{fig_second_caseA3}}} 
\caption{(a) Aggregated demand of price-responsive users before reshaping, $ P_\text{U}^\text{NSW}$, aggregated flexible demand, $ P_\text{F}^\text{NSW}$, and aggregated PV generation, $ P_\text{PV}^\text{NSW} $,  (b) battery charging, discharging and SOC in NSW during 30\superscript{th} and 31\superscript{st} of January 2020.}
\label{figure:Efficacy}
\end{figure}

 Fig.~\ref{figure:Efficacy} demonstrates the results of the demand model in NSW for low DR Scenario during the 30\superscript{th} and 31\superscript{st} of January 2020. As it can be seen in Fig.~\ref{figure:Efficacy}a, flexible demand, $ P_\text{F}^\text{NSW} $, increases during cheap hours in comparison with $ P_\text{U}^\text{NSW} $ because price-responsive users use their battery storage to store cheaper electricity (Fig.~\ref{figure:Efficacy}b). However, flexible demand decreases during peak hours due to PV generation, $ P_\text{PV}^\text{NSW} $, and also battery storage discharge (Fig.~\ref{figure:Efficacy}b).  Also, Fig.~\ref{figure:Efficacy}a clearly illustrates that $\sum_{h=1}^{H}P_\text{F}^m(h)\Delta h= \sum_{h=1}^{H}(P_\text{U}^m(h)-P_\text{PV}^m(h))\Delta h+E_\text{loss}^{\text{B},m}$. Furthermore, Fig.~\ref{figure:Efficacy}b shows that with appropriate selection of the $P_\text{B,cha}^{\text{max},m}$, battery SOC remains within the bounds. This shows the efficacy of the proposed model.It should be noted that all the results in this section assume that price anticipating users do not sent power back to the grid ($P_\text{L}^{m}\leq P_\text{LF}^{m}$).


\subsection{Comparison of the load profile with different approaches}

\par
\begin{figure}
{\subfloat[]{\includegraphics[width=8.8cm, height=6.5cm]{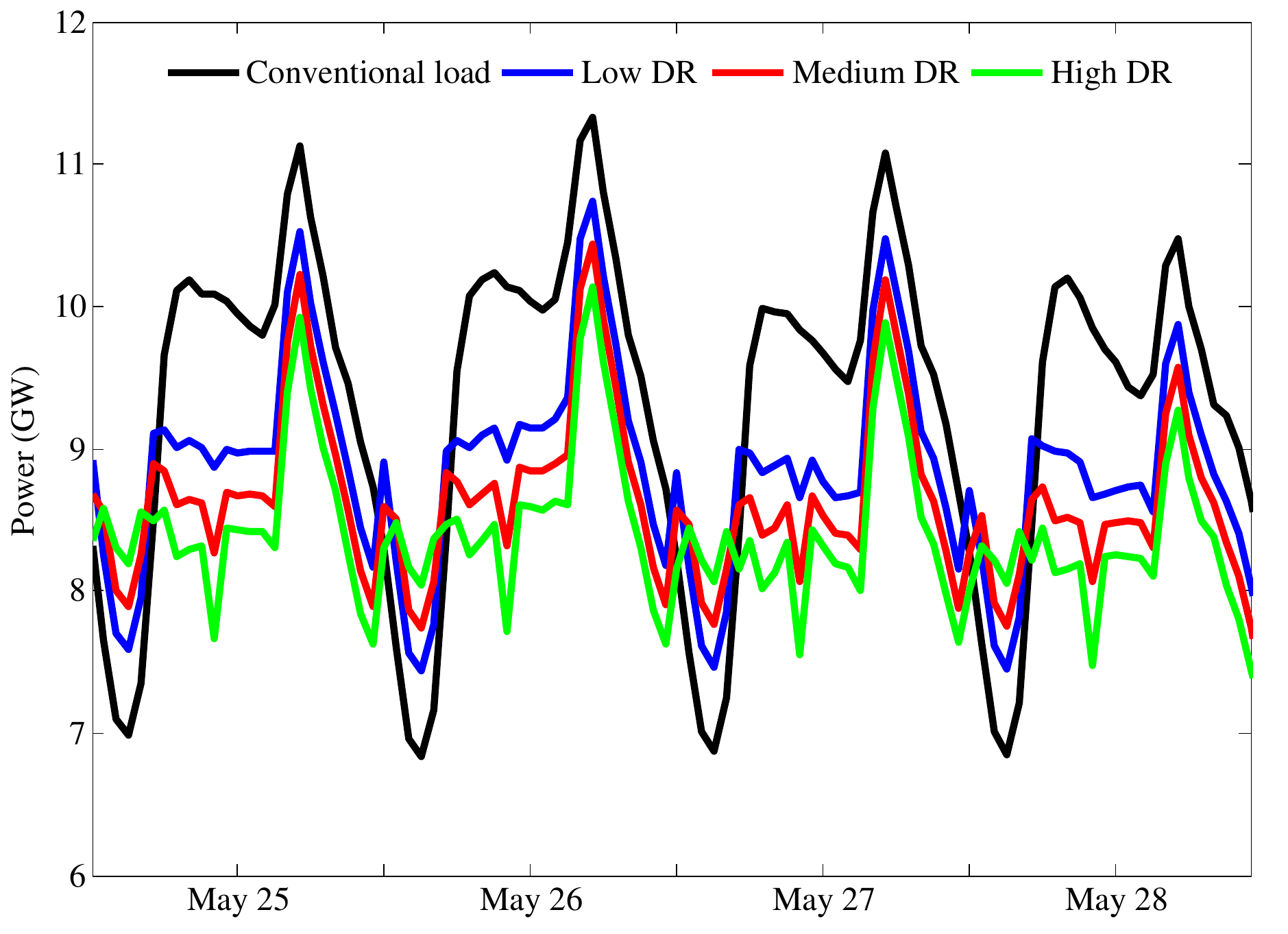}
\label{fig:first_caseA111}}}
{\subfloat[]{\includegraphics[width=8.8cm, height=6.5cm]{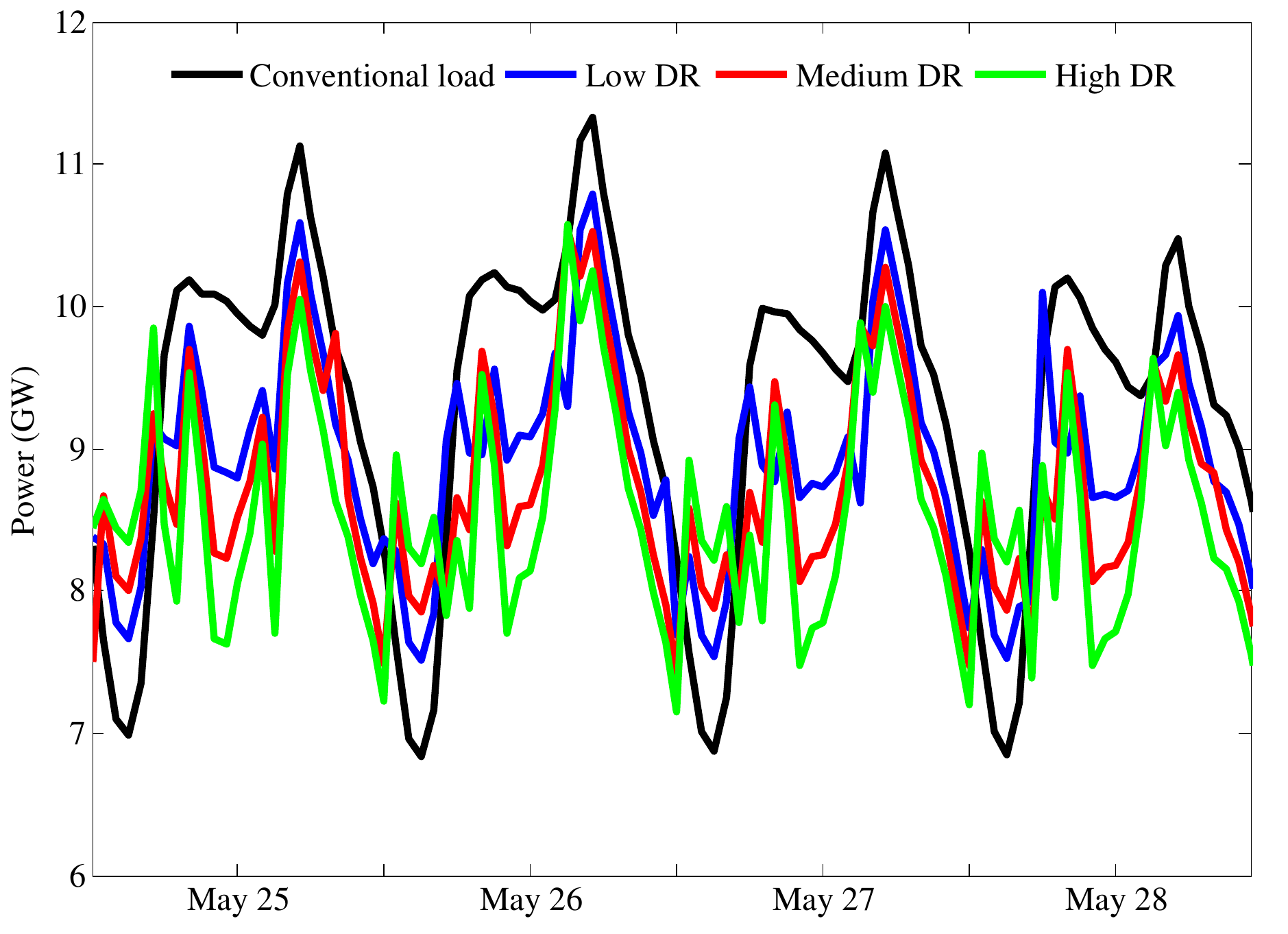}
\label{fig_second_caseA112}}} 
\caption{Demand profile in NSW during 25\superscript{th} to 28\superscript{th} of May 2020 with (a) price anticipating and (b) price taking assumptions.}
\label{figure:The effect of different DR scenarios on the load profile}
\end{figure}

 Figs.~\ref{figure:The effect of different DR scenarios on the load profile}a and \ref{figure:The effect of different DR scenarios on the load profile}b compare the effect of price anticipating and price taking assumptions, and also different levels of DR on the load profile resulting from solving the proposed model and the demand model in \cite{Hesam2014} in NSW during the 25\superscript{th} to 28\superscript{th} of May 2020, respectively. As it is shown in the Figs.~\ref{figure:The effect of different DR scenarios on the load profile}a and \ref{figure:The effect of different DR scenarios on the load profile}b, in both approaches, the users shift their consumption (using PV-storage system) from expensive hours to cheaper ones to utilise zero cost electricity produced by RESs. In the proposed model, the effect of the users' action on the electricity price signal is considered by the load aggregators. As a result, the load profile with price anticipating assumption (Fig.~\ref{figure:The effect of different DR scenarios on the load profile}a) is smoother than the load profile with price taking assumption (Fig.~\ref{figure:The effect of different DR scenarios on the load profile}b). In the latter approach the effect of users' action is not considered on the electricity price signal. So, price-responsive users may shift their consumption to cheaper time slots all-together, which results in secondary load peaks due to load synchronisation, as it can be seen in Fig.~\ref{figure:The effect of different DR scenarios on the load profile}b. 

\par In the next subsection, the effect of the both demand models on balancing and loadability of the NEM in 2020 are studied with the increased penetration of RESs.

\subsection{Balancing and loadability results}

\subsubsection{BAU Scenario}

\par In the BAU Scenario, a big portion of the demand is supplied by coal-fired power plants and the peak loads are met with backup supply (i.e. gas turbines (GTs)). The results for spilled energy and hours, supplied electrical energy from GTs and average loadability for all the scenarios over the simulated year are summarised in Table~\ref{tab:Balancing, loadability and minimum singular value  of $G_{S}$ for scenarios 1 to 4}. For the loadability calculation, we checked different load/generation increase scenarios. Among them, we showed the results for one of the critical ones in Table~\ref{tab:Balancing, loadability and minimum singular value  of $G_{S}$ for scenarios 1 to 4} in which all loads in QLD are assumed to increase uniformly in small steps with constant power factor. Also, it is assumed that all the generators in QLD are scheduled with the same participation factor to pick up the system loads. The loadability is computed for each hour until a step before power flow divergence.

\subsubsection{Renewable Scenarios}

\par In the rest of the paper, Renewable Scenarios with the conventional load, the load model in \cite{Hesam2014} and the proposed model are called CL, PTDR and PADR, respectively. Also, PTDR and PADR scenarios with low, medium and high DR penetrations are called PTDR1-PTDR3 and PADR1-PADR3 respectively in Table~\ref{tab:Balancing, loadability and minimum singular value  of $G_{S}$ for scenarios 1 to 4}. Unserved hours for all scenarios are zero. Comparing the BAU Scenario and the Renewable Scenario with the conventional load (CL), it can be seen that with the increased penetration of RESs, the loadability is reduced from \SI{27.13}{\giga\watt} to \SI{22.17}{\giga\watt}.Also, the required electrical energy from the backup generation (i.e. GTs) is increased from \SI{18.73}{\tera\watt}to \SI{18.77}{\tera\watt}.

\par Compared to the CL Scenario (i.e. zero DR penetration), higher penetrations of DR (PADR and PTDR) improve the performance and loadability, and reduce the required energy from backup supply. High DR Scenario with price anticipating assumption (PADR3) has the highest loadability followed by the medium and low DR Scenarios. The loadability is increased from \SI{22.17}{\giga\watt} for the CL Scenario to \SI{26.24}{\giga\watt} for the PADR3 Scenario, which implies a considerable improvement in loadability. However, with price taking assumption, medium DR Scenario (PTDR2) has the highest loadability (i.e. \SI{25.53}{\giga\watt}). Increasing the penetration of DR with price taking assumption beyond a certain point (from PTDR2 to PTDR3), fails to improve loadability further. This is because of secondary peaks for high DR scenario which deteriorates loadability compared to lower DR penetrations. Furthermore, in both approaches, high DR Scenario has the lowest spilled energy and hours followed by medium and low DR Scenarios. Also, the required energy from the backup generation for DR scenarios with price taking assumption is slightly higher than DR scenarios with price anticipating assumption. This happens as the price taker users may shift their consumption to the cheaper time slots all-together, and therefore some backup supply might be required to keep the network in balance.

\begin{table}
\centering
\caption{balancing and loadability results for all the scenarios}
\begin{tabular}{| c | c | c | c | c |}\hline
Scenarios & Spilled energy  & Spilled hours  & GT energy & Loadability \\
      &   (TWh)    &(\%)      &(TWh)& (GW) \\\hline
BAU   &  -   &   -   & 18.73 & 27.13  \\
CL    & 0.71 & 13.65 & 18.77 & 22.17  \\
PADR1 & 0.66 & 13.03 & 18.28 & 24.82  \\
PADR2 & 0.61 & 12.67 & 17.78 & 25.86  \\
PADR3 & 0.54 & 11.71 & 17.04 & 26.24  \\
PTDR1 & 0.66 & 13.03 & 18.33 & 23.78  \\
PTDR2 & 0.61 & 12.67 & 17.85 & 25.53  \\
PTDR3 & 0.54 & 11.71 & 17.28 & 24.41  \\\hline
\end{tabular}
\label{tab:Balancing, loadability and minimum singular value  of $G_{S}$ for scenarios 1 to 4}
\end{table}

\par
\begin{figure} 
{\subfloat[]{\includegraphics[width=8.8cm, height=6.6cm]{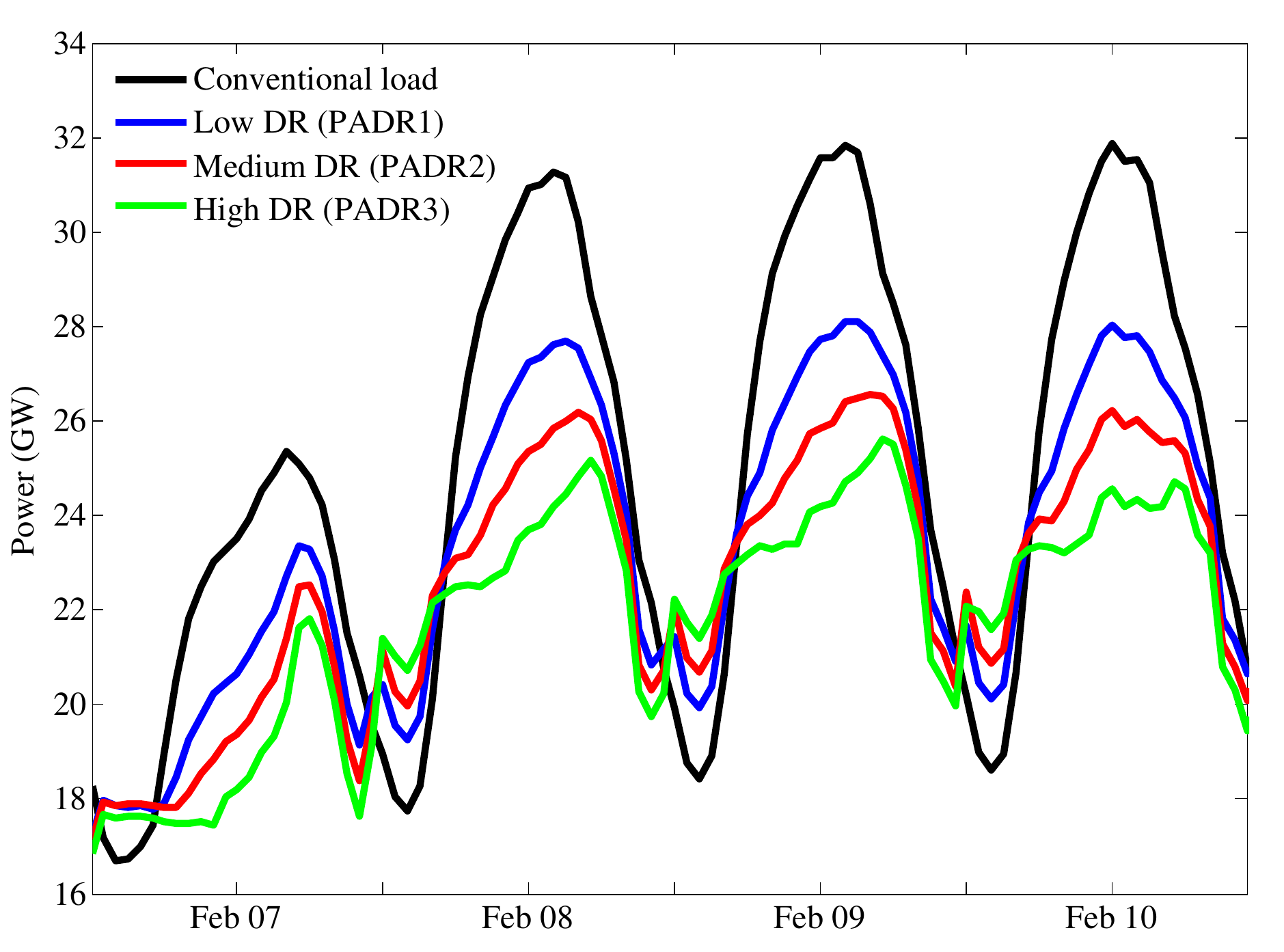}
\label{fig:first_caseA211}}}
{\subfloat[]{\includegraphics[width=8.8cm, height=6.5cm]{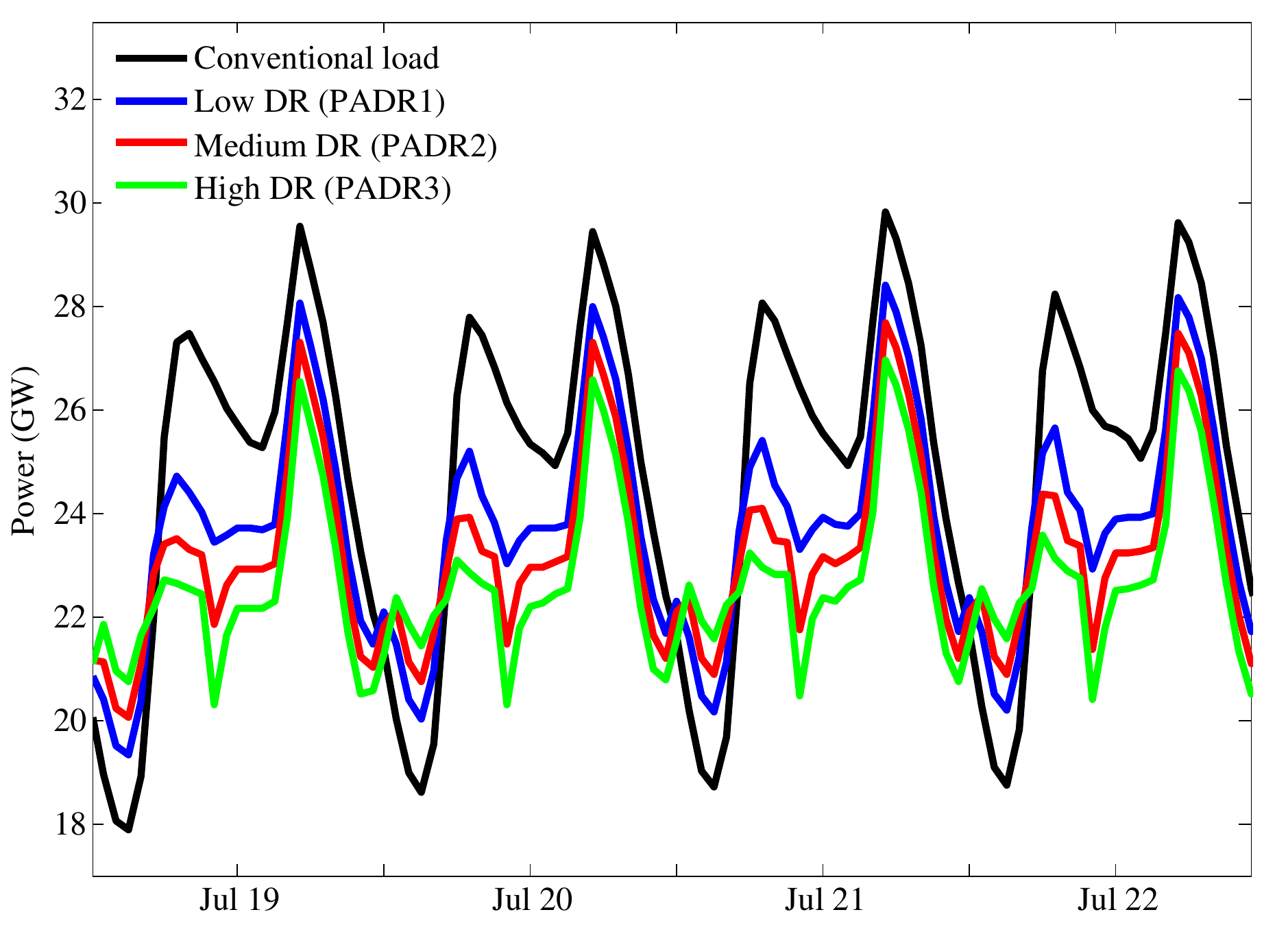}
\label{fig_second_caseA212}}} 
\caption{Demand profile for the NEM for Scenarios CL and PADR1-3 in one of the critical (a) summer and (b) winter peaks.}
\label{figure:Demand profile for scenarios 2 to 5 in one of the critical (a) summer and (b) winter peak}
\end{figure}

\begin{figure} 
{\subfloat[]{\includegraphics[width=8.8cm, height=6.6cm]{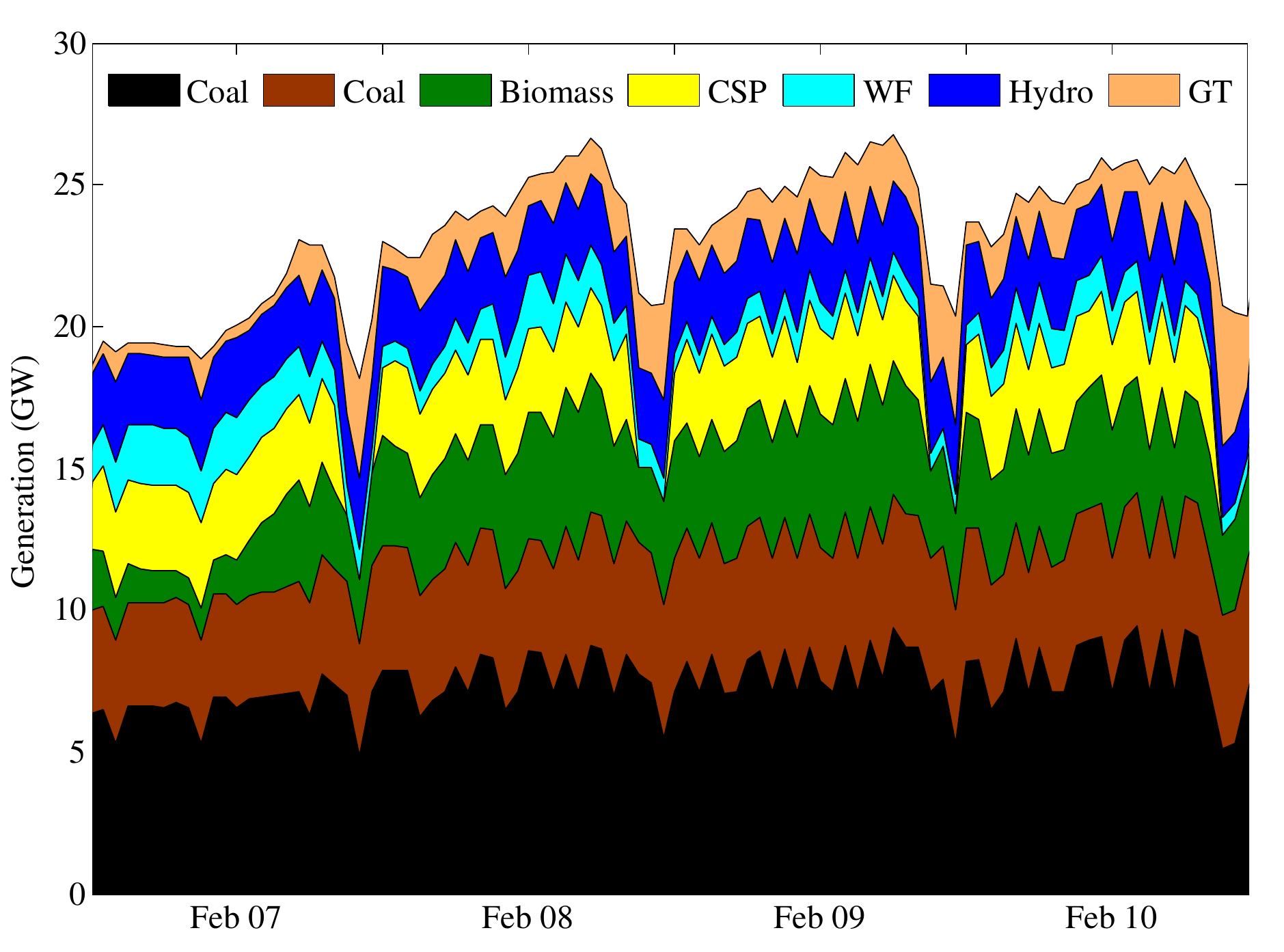}
\label{fig:first_caseA311}}}
{\subfloat[]{\includegraphics[width=8.8cm, height=6.5cm]{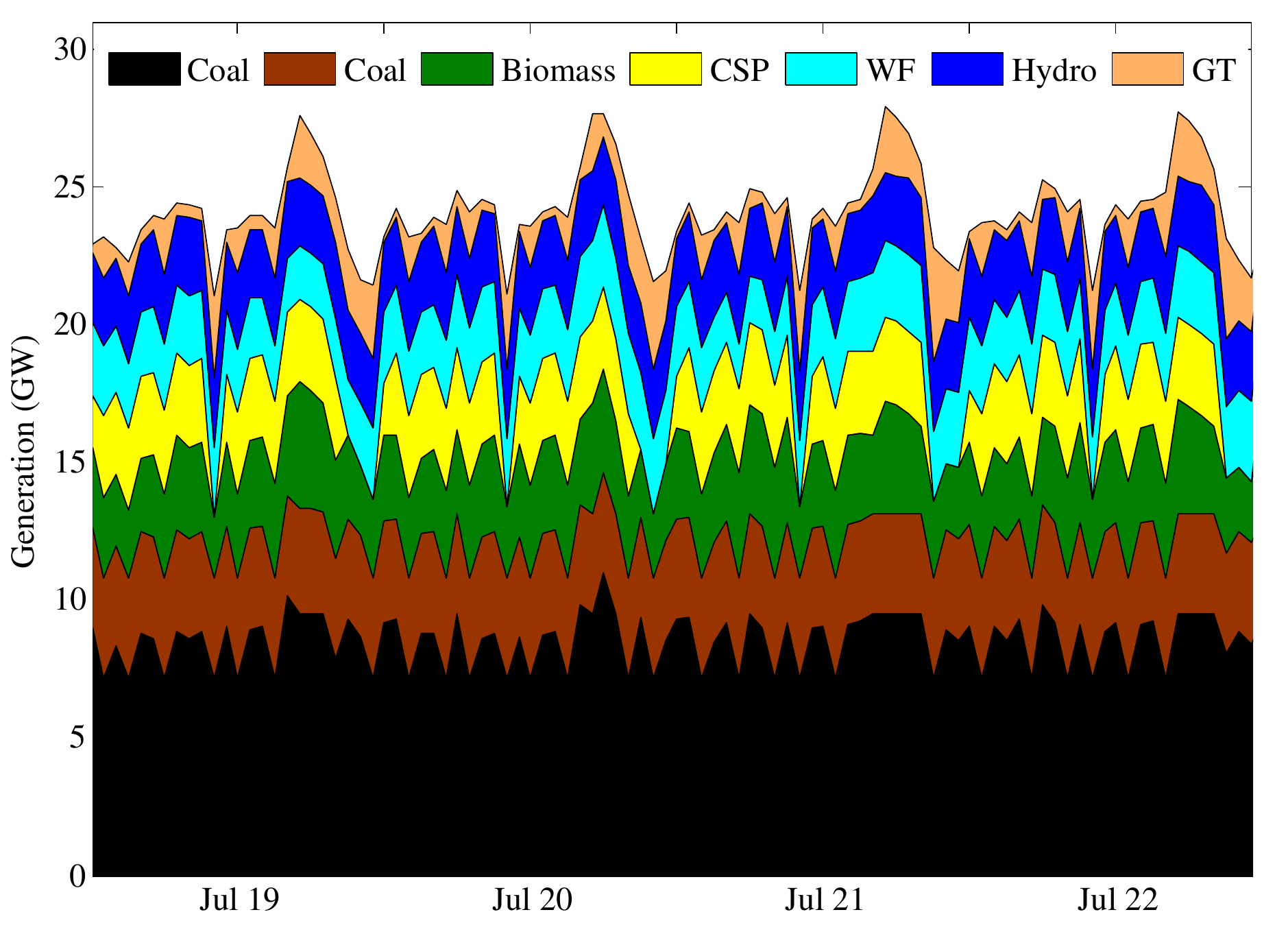}
\label{fig_second_case312}}}
\caption{Balancing results for the NEM for the high DR Scenario with price anticipating assumption (PADR3) in one of the critical (a) summer and (b) winter peaks.}
\label{figure:Balancing results for scenario 3 from 10th to 15th of January 2020}
\end{figure}

 Fig.~\ref{figure:Demand profile for scenarios 2 to 5 in one of the critical (a) summer and (b) winter peak} shows the the effect of different DR scenarios on the load profile when users are price anticipators in one of the critical summer (i.e. 07\superscript{th} to 10\superscript{th} of February) and winter peaks (i.e. 19\superscript{th} to 22\superscript{th} of July) for the NEM. The balancing results for the high DR scenario during those peak hours are demonstrated in Fig.~\ref{figure:Balancing results for scenario 3 from 10th to 15th of January 2020}. As shown in Fig.~\ref{figure:Balancing results for scenario 3 from 10th to 15th of January 2020}, in the critical summer days the wind is not strong enough and the output of WF is low, and in critical winter days the solar exposure reduces and CSP output decreases as well. However, due to enough response capacity from DR (Fig.~\ref{figure:Demand profile for scenarios 2 to 5 in one of the critical (a) summer and (b) winter peak}) and backup supply, balancing and loadability under these worst-case conditions are maintained. Figs.~\ref{figure:Demand profile for scenarios 2 to 5 in one of the critical (a) summer and (b) winter peak} and~\ref{figure:Balancing results for scenario 3 from 10th to 15th of January 2020} show that DR can help balance fluctuating RESs and demand in real-time even under the worst-case conditions.

\section{Conclusion}

\par This paper proposes a demand model which integrates the aggregated effect of DR in a simplified representation of the effect of market/dispatch processes. The model is intended to be used for system studies at transmission levels in which users are assumed to be price anticipators. In the proposed optimisation formulation, the conventional demand model is augmented by including the aggregated effect of numerous price anticipating users equipped with PV-storage systems. As a case study, the effect of the demand model on performance and loadability of the NEM in 2020 with the increased penetration of RESs is studied using a 14-generator model. Also, the results are compared with the demand model in which users are assumed to be price takers.

\par Simulation results show that with the increased penetration of RESs and no DR loadability is reduced. With DR, however, loadability is improved and the required backup supply is reduced. Increasing the penetration of DR with price anticipating assumption improves loadability. However, when users are price takers, increasing the penetration of DR beyond a certain point does not necessarily improve the loadability and might even deteriorate it. This is due to load synchronisation of price taking users, which might results in secondary peaks.

\par The future research aims at expanding the idea in this paper by using bi-level optimisation to avoid implicit modelling of battery storage and heuristic search. Also, the effect of the aggregated demand model including DR will be studied on performance and stability of FG scenarios.

\ifCLASSOPTIONcaptionsoff
  \newpage
\fi



%



\bibliographystyle{IEEEtran}
\bibliography{Research_Papers}

%








\end{document}